\documentclass[12pt,a4paper]{amsart}
\newif\ifpdf 
\ifx\pdfoutput\undefined 
        \pdffalse 
\else 
        \pdfoutput=1 
        \pdftrue 
\fi

\ifpdf 
        \usepackage[pdftex]{graphicx} 
        \DeclareGraphicsExtensions{.pdf}
        
\else 
        \usepackage{graphicx}
        \DeclareGraphicsExtensions{.eps}
\fi
\usepackage{mathrsfs}
\usepackage{bbm}
\usepackage{enumerate}
\usepackage{moreverb}

\newcommand\pref[1]{(\protect\ref{#1})}
\newcommand\bdm{\begin{displaymath}}
\newcommand\edm{\end{displaymath}}

\newcommand\ba{{\mathbb A}}
\newcommand\be{{\mathbb E}}
\newcommand\br{{\mathbb R}}
\newcommand\bc{{\mathbb C}}
\newcommand\bq{{\mathbb Q}}
\newcommand\bp{{\mathbb P}}

\newcommand\bz{{\mathbb Z}}
\newcommand\bbb{{\mathbb B}}
\newcommand\by{{\mathbf y}}
\newcommand\bff{{\mathbb F}}

\newcommand\cc{{\mathcal C}}
\newcommand\ch{{\mathcal H}}
\newcommand\fd{{\mathfrak D}}
\newcommand\mm{{\mathscr M}}
\newcommand\nn{{\mathscr N}}
\newcommand\cals{{\mathscr S}}
\newcommand\gog{{\mathfrak g}}

\DeclareMathOperator{\aut}{Aut}

\newtheorem{thm}{Theorem}[section]
\newtheorem{zvk1}[thm]{Theorem of Zariski-Van Kampen}
\newtheorem{prop}[thm]{Proposition}
\newtheorem{cor}[thm]{Corollary}
\newtheorem{lema}[thm]{Lemma}
\newtheorem{fact}[thm]{Fact}
\newtheorem{thm0}{Theorem}
\theoremstyle{remark}
\newtheorem{obs}[thm]{Remark}
\theoremstyle{definition}
\newtheorem{dfn}[thm]{Definition}
\newtheorem{ntc}[thm]{Notation}

\newtheorem{ejm}[thm]{Example}
\newtheorem{str}[thm]{Strategy}
\newtheorem{cnt}[thm]{Construction}
\numberwithin{equation}{section}

\title[Braid monodromy and plane curves]
{Effective invariants of braid monodromy
and topology of plane curves}

\author[E. Artal]
{Enrique ARTAL BARTOLO}
\address{Departamento de Matem\'aticas\\
Campus Plaza de San Francisco s/n\\
E-5009 Zaragoza SPAIN}
\email{artal@posta.unizar.es}

\author[J. Carmona]{Jorge CARMONA RUBER}
\address
{Departamento de Sistemas inform\'aticos y programaci\'on\\
Universidad Complu\-tense\\
Ciudad Universitaria s/n\\
E-28040 Madrid SPAIN}

\email{jcarmona@eucmos.sim.ucm.es}

\author[J.I. Cogolludo]{Jos\'e Ignacio COGOLLUDO AGUST\'IN}
\address{Departamento de \'Algebra\\
Universidad Complutense\\
Ciudad Universitaria s/n\\
E-28040 Madrid SPAIN}

\email{jicogo@eucmos.sim.ucm.es}
\date\today

\thanks{First and second authors are partially supported by
DGES PB97-0284-C02-02; third author is partially
supported by DGES PB97-0284-C02-01}

\begin{document}
\begin{abstract} In this paper we construct 
effective invariants for braid monodromy of affine curves.
We also prove that, for some curves, braid monodromy
determines their topology. We apply this result
to find a pair of curves with conjugate equations
in a number field but which do not admit
any orientation-preserving
homeomorphism.
\end{abstract}

\maketitle

Let $\cc \subset \bc^2$ be an algebraic affine curve.
We say a property $P(\cc)$ is an invariant of $\cc$ if it 
is a topological invariant of the pair $(\bc^2,\cc)$, in
other words, if $P(\cc)=P(\cc')$ whenever $(\bc^2,\cc)$ 
and $(\bc^2,\cc')$ are homeomorphic as pairs. Analogously,
we define the concept of invariants of projective algebraic
curves.

Our purpose in this paper is the construction of new
and effective invariants for algebraic curves that
reveal that the position of singularities is not enough
to determine the topological type of the pair 
$(\bp^2,\cc)$.

These invariants will be derived from a refinement of a
well-known invariant of curves such as braid monodromy.
Roughly speaking, the braid monodromy of a curve 
$\tilde \cc \subset \bp^2$ with respect to a pencil of 
lines $\ch$, is a representation of a free group 
$\bff$ on the braid group on $d$ strings --\,where $d$ 
is the degree of $\tilde \cc$ restricted to the generic 
fiber of $\ch$ after resolution of its base point.
The free group $\bff$ corresponds to the fundamental 
group of an $r$-punctured complex line, where the punctures 
come from the non-generic elements of $\ch$.

Braid monodromy is a strong invariant for projective
(or affine) plane curves. It is fair to say that the main
ideas leading to this invariant have already appeared in the 
classic works of Zariski \cite{zr:29} and Van Kampen 
\cite{vk:33} as a tool to find the fundamental group of 
the complement of a curve.
Nevertheless, its consideration as an invariant itself 
is due to B.~Moishezon --\,see \cite{mz:81} for definitions 
and other related references. This author defines the
invariant and obtains from it a number of beautiful results 
in several papers, both as single author and with M.~Teicher.
Later on, A.~Libgober \cite{li:86} proves that the homotopy 
type of the complement of an affine curve is an invariant 
of its braid monodromy. Several invariants of the
fundamental group have been found to be effective as 
invariants of affine curves, such as the Alexander
polynomial, the Alexander module or the sequence of 
characteristic varieties \cite{li:98}. A second group of 
effective invariants was described by A.~Libgober in 
\cite{li:89}. They depend on polynomial representations of 
the braid group and are invariants of the conjugation class 
of the image of the representation defining the braid monodromy. 
\smallbreak

By fixing a particular class of bases of the aforementioned 
free group $\bff$ --\,a \emph{geometric basis}\,-- braid 
monodromy may be represented by an $r$-tuple of braids.
This sequence is not unique; there is a natural action
of $B_d\times B_r$ on $B_d^r$ such that the admissible
$r$-tuples form an orbit of this action. In general, it is 
not easy to decide whether two elements of $B_d^r$ are in 
the same orbit or not, since these orbits are infinite.
\smallbreak

Our purpose is to find finer invariants that are sensitive 
to changes within conjugation classes. Using finite 
representations of braid groups we obtain new effective 
invariants for braid monodromies. 
\emph{Effectiveness} is obtained 
by means of the free software GAP4 \cite{GAP4}.

\smallbreak
For curves having only nodes and cusps as singularities,
Teicher and Kulikov \cite{kt:00} have recently proved
that the diffeomorphism type of their embedding
is determined by their braid monodromy.
The general case has been resolved independently by the
second named author in \cite{car:xx}. In this work, we 
prove a kind of converse of this result; we will prove 
that braid monodromy of an affine curve $\cc$
determines the oriented topological type of the
pair $(\bp^2,\overline{\cc^\varphi}\cup L_\infty)$, 
where $\cc^\varphi$ is the
union of $\cc$ and all non-transversal vertical lines
and $L_\infty$ is the line at infinity.

\smallbreak
Braid monodromy may also be useful to study the moduli
of curves with prescribed degree and topological types
of singularities. Although braid monodromy is defined for 
affine curves, one can define it in the projective case
by choosing generic lines at infinity. Let $d$ be a 
positive integer and let $T_1,\dots,T_r$ be topological 
types of singularities of curves. Let us denote by 
$\Sigma(T_1,\dots,T_r;d)$ the space of all projective 
plane curves of degree $d$ with $r$ singular points of 
type $T_1,\dots,T_r$. Let $\mm(T_1,\dots,T_r;d)$ be the 
quotient of $\Sigma(T_1,\dots,T_r;d)$ by the action
of the projective group. Sometimes we will restrict 
ourselves to $\Sigma^{\text{irr}}$ and $\mm^{\text{irr}}$
by considering only irreducible curves.
Note that braid monodromy is an invariant
of each connected component of these moduli spaces.

\smallbreak
In this work we also address the problem of the topology
type of conjugate varieties. Let $(V,W)$ be a pair
of projective varieties, $W\subset V$ ($W$ may be empty),
with defining equations in a number field $K$.
For each embedding $j$ of $K\subset\bc$ (if we consider $K$
already embedded in $\bc$, for each Galois action
of the normal closure of $K/\bq$), one obtains
a pair $(V_j,W_j)$ of complex projective varieties,
which essentially shares all the algebraic properties
of $(V,W)$.
\medbreak

It is well known that a great number of topological
properties of $(V_j,W_j)$ are of algebraic nature
and they depend only on the pair $(V,W)$.
Nevertheless, examples by Serre \cite{se:64} and 
Abelson \cite{ab:74} show that it is possible to find 
examples of non-homeomorphic conjugate varieties, i.e.,
different embeddings of $K$ provide different topological 
pairs.

\medbreak
Serre distinguishes the complex manifolds
by means of the fundamental group (although
the associated algebraic fundamental groups
must be isomorphic). Abelson's examples
have the same fundamental groups but
they differ on the cup product in cohomology.
Serre's examples are surfaces whereas Abelson's
are in higher dimensions. Applying
generic projection and Chisini's conjecture,
Serre's result implies that there are conjugate projective
plane curves which do not have
the same embedded topological type in the projective plane,
probably with a large degree.
Recently, A.~Kharlamov and V.~Kulikov \cite{khku:01}
have proven the existence of pairs of (complex) conjugate
algebraic curves which are of course diffeomorphic but
not isotopic. Note that the homeomorphism does not
respect curve orientations. In addition, they prove
that their examples have non-equivalent
braid monodromies; the smallest degree of their
examples is 825.

\medbreak
We want to construct such examples with smaller degrees.
It is well known, from Degtyarev's work \cite{deg:90}, that
no example of non-homeomorphically embedded conjugate 
curves can exist up to degree $5$. In \cite{acc:00}, 
several examples were found as candidates to produce such 
kind of curves.
For example, it was proved that $\mm(\ba_{19};6)$ has 
two elements and representatives were found having 
conjugate equations in $\bq(\sqrt{5})$. In other
cases, such as $\mm(\ba_{18},\ba_1;6)$, there are three 
elements conjugated in a degree $3$ extension of $\bq$.
In \cite{acct:00}, we have also studied the family
$\mm(\ba_9,\ba_9,\ba_1;6)$. It is easily seen that if 
$\cc\in\mm(\ba_9,\ba_9,\ba_1;6)$ then $\cc$ is the union 
of a quintic curve with two singular points 
(one $\ba_9$ and one $\ba_{1}$) and a line which
intersects the quintic at a smooth point
with intersection multiplicity $5$
(providing the \emph{second} $\ba_9$ point).
Following the ideas used in \cite{acct:00}
we can prove that $\mm(\ba_9,\ba_9,\ba_1;6)$
has three points; the first one contains curves
such that  the tangent line to $\cc$ at
the first $\ba_9$ point passes through the
second $\ba_9$ point. This is not the case for
the last two cases, where we  can find representatives
having conjugate equations with coefficients in $\bq(\sqrt{5})$.
H. Tokunaga has proven that curves in the first case
are not topologically equivalent
to those in the other ones, by means of
finding in an algebraic way all possible coverings
corresponding to the dihedral group of $10$ elements.

In this work we will deal with curves in
$\mm:=\mm(\be_6,\ba_7,\ba_3,\ba_2,\ba_1;6)$. It is easily
seen that if $\cc\in\mm$ then $\cc$ is the union of
a quintic curve with three singular points (of types
$\be_6$, $\ba_3$ and $\ba_{2}$) and a line which
intersects the quintic at two smooth points
with intersection multiplicities $4$ and $1$.

\begin{thm0} The space $\mm$ has exactly two points $\mm_1$,
and $\mm_2$. There exist representatives $\cc_i\in\mm_i$, $i=1,2$,
having conjugate equations with coefficients in $\bq(\sqrt{2})$.
\end{thm0}

We will not give a complete proof of this theorem
since it is obtained by applying a standard
Cremona transformation and the ideas used in 
\cite{acct:00} to the family $\mm(\ba_{15},\ba_3,\ba_1;6)$. 
In this paper we also compute the
necessary tools to prove that the fundamental groups of 
the complements are both isomorphic to $\bz\times SL(2,\bz/7\bz)$.

\medbreak
Let $\cc\in\mm$; we construct a curve $\tilde \cc$
of degree $12$ as follows: let us consider the pencil
of curves through the $\be_6$ point. There are exactly
six non-generic lines in this pencil which are: 
the lines joining this point with the four other singular 
points of $\cc$, an ordinary tangent line and the tangent 
line to $\cc$ at the base point. Then, $\tilde \cc$ is the union 
of the quintic and seven lines. Note that if 
$\cc_i\in\mm_i$, $i=1,2$, are conjugate curves, it is also the 
case for $\tilde \cc_i$. We will prove the following:

\begin{thm0}\label{conj} Let us consider two conjugate curves
$\cc_i\in\mm_i$, $i=1,2$. Then, $\tilde \cc_1$ and $\tilde \cc_2$
are two conjugate curves such that the pairs
$(\bp^2,\tilde \cc_1)$ and $(\bp^2,\tilde \cc_2)$
are not homeomorphic by an orientation-preserving
homeomorphism.
\end{thm0}

In order to prove this theorem, we construct
braid monodromies for the affine curves resulting from
$\cc\in\mm_1\cup\mm_2$ and considering the tangent line
at the $\be_6$ point as the line at infinity.
The pencil used for this purpose is the one given by
the lines through the $\be_6$ point. Braid monodromies 
of curves in different connected components turn out 
not to be equivalent --\,note that this is also the case
in the examples by Serre and Kharlamov-Kulikov.
\smallbreak

In \S 1, we set notations and definitions.
In \S 2, we define the concepts of geometric and
lexicographic bases and lexicographic braids. 
In \S 3, we translate the definition of braid monodromy 
in terms of $r$-tuples of braids and prove
the theorem about the determination
of the braid monodromy by the topology
in the fibered case.
In \S 4 and \S 5, we introduce a way to produce
effective invariants for braid monodromies and
finite representations of braid groups.
In \S 6, we find braid monodromies for representatives 
in $\mm_1$ and $\mm_2$. In the Appendix, we give
the GAP4 programs which help to distinguish braid
monodromies. The results of \S 3, \S 6 and
the Appendix prove Theorem \ref{conj}.

\section{Settings and definitions}

Let $G$ be a group and $a,b\in G$. For the sake of 
simplicity we use the following notation
$a^b:=b^{-1}ab$, $[a,b]:=a^{-1}b^{-1}ab$ and
$b* a:=bab^{-1}$.
Given any set $A$ we denote by $\Sigma_A$ the group 
of all bijective maps from $A$ onto itself. As usual
we denote $\Sigma_{\{1,\dots,n\}}$ by $\Sigma_n$.

\smallbreak
We also denote by $\bff_n$ the free group on $n$ generators, 
say $x_1,\dots,x_n$, and by $B_n$ the braid group on $n$ 
strings given by the following presentation:
\bdm
\left\langle \sigma_1,\dots,\sigma_{n-1} :\
[\sigma_i,\sigma_j]=1,\ |i-j|\geq 2,
\sigma_i\sigma_{i+1}\sigma_i=\sigma_{i+1}\sigma_i\sigma_{i+1},\ 
i=1,\dots,n-2
\right\rangle.
\edm
There is a natural right action
$\Phi:\bff_n\times B_n\to\bff_n$ such that
\bdm
x_j^{\sigma_i}=
\begin{cases}
x_{i+1}&\text{if $j=i$}\\
x_{i+1}* x_i&\text{if $j=i+1$}\\
x_j&\text{if $j\neq i,i+1$}.
\end{cases}
\edm

The induced antihomomorphism $B_n\to\aut\bff_n$ is injective.
Identifying $B_n$ with its image one has the following
characterization, which may be found in any classic
text on braids, e.g. 
%\cite{art:26}, 
\cite{art:47} or \cite{bir:74}:

\begin{prop} Let $\tau:\bff_n\to\bff_n$ be an automorphism. 
Then $\tau\in B_n$ if and only if 
$\tau(x_n\cdot\ldots\cdot x_1)=x_n\cdot\ldots\cdot x_1$
and there exists $\sigma\in\Sigma_n$ such that
$\tau(x_i)$ is conjugate to $x_{\sigma(i)}$.
\end{prop}

Next we describe a geometrical interpretation
of these definitions. All the results
in this section appear in J.~Birman's book \cite{bir:74}. 
Let us fix 
$\by:=\{t_1,\dots,t_n\}$ a subset of $\bc$ of exactly 
$n$ elements. Let us consider also a \emph{big enough} 
geometric closed disk $\Delta$ such that $\by$ is 
contained in the interior of $\Delta$ and let $*$ be
a point on $\partial\Delta$. 

\begin{ntc}\label{vuelta} We will denote by
$\gamma\in\pi_1(\bc\setminus\by; *)$
the homotopy class of the loop based at $*$
which surrounds the circle $\partial\Delta$
counterclockwise.
\end{ntc}

We recall:

\begin{fact} The group $\pi_1(\bc\setminus\by; *)$ is isomorphic
to $\bff_n$.
\end{fact}

\begin{dfn} A \emph{$\by$-special homeomorphism} is an 
orientation-preserving homeomorphism $\bc\to\bc$ which globally 
fixes $\by$ and is the identity on $\overline{\bc\setminus\Delta}$.
A \emph{$\by$-special isotopy} 
is an isotopy $H:\bc\times[0,1]\to\bc\times[0,1]$
with the above properties.
\end{dfn} 

The quotient of the
set of $\by$-special homeomorphisms
module $\by$-special isotopy has a natural
structure of group and is denoted by $\bbb_{\by}$.
There is a natural left action of $\bbb_{\by}$ on
$\pi_1(\bc\setminus\by; *)$ as a result of the definitions.

\begin{fact} The group $\bbb_{\by}$ is isomorphic
to $B_n$.
\end{fact}

Let $V$ be the space of monic polynomials in $\bc[t]$
of degree $n$ ($V$ and $\bc^n$ are naturally isomorphic as
affine spaces). Let $D$ be the discriminant
space of $V$, i.e., the set of elements in $V$ with
multiple roots. The space $D$ is an algebraic hypersurface of $V$.
We will use the following notation: $X:=V\setminus D$, 
$p(t):=\prod_{j=1}^n(t-t_j)$, $B_{\by}:=\pi_1(X;p(t))$. 
By means of the theorem of continuity of roots, $X$
can be naturally identified
with the subsets of $\bc$ having exactly $n$ elements.

\begin{obs}
From now on, $X$ will denote, depending on the context, 
either $V\setminus D$ or the set 
$\{\by\subset\bc\ |\ \#\by=n\}$ with quotient topology
resulting from the $n^{\text{th}}$-symmetric product of $\bc$.
\end{obs}
\smallbreak

An element $\tau\in B_{\by}$ is identified with
the homotopy class, relative to $\by$, of sets of paths 
$\{\gamma_1,\dots,\gamma_n\}$, $\gamma_j:[0,1]\to\bc$, 
starting and ending at $\by$ and such that
$\forall t\in[0,1]$, $\{\gamma_1(t),\dots,\gamma_n(t)\}$
is a set of $n$ distinct points. The elements of
$B_{\by}$ are called braids based at $\by$ and
are represented as a set of non-intersecting paths in 
$\bc\times[0,1]$, as usual. 
Products are also defined in the standard way.
\smallbreak

Since any orientation-preserving homeomorphism
of $\bc$ is isotopic to the identity,
we can associate to any element of $\bbb_{\by}$
the braid representing the motion of $\by$ along
the isotopy.

\begin{fact} The groups $\bbb_{\by}$ and
$B_{\by}$ are naturally anti-isomorphic.
\end{fact}

Thus one can construct a right action
\bdm
\Phi_{\by}:\pi_1(\bc\setminus\by; *)\times B_{\by}\to
\pi_1(\bc\setminus\by; *).
\edm
As above, exponential notation will be used to describe
this action.

\begin{ejm}\label{standard} Let us consider:
\begin{itemize} 

\item $\by^0=\{-1,\dots,-n\}$,

\item $\Delta=\{t\in\bc|\ |t|\leq R\}$, $R\gg 0$,

\item $*=R$.
\end{itemize}
 
Let $\delta_j$ be the path obtained
by surrounding counterclockwise the circle of radius 
$\frac{1}{4}$ centered at $-j$ and starting at 
$-j+\frac{1}{4}$, $j=1,\dots,n$.
The path $\delta_j^+$ will denote the first half-circle 
and $\delta_j^-$ the second one.

For any given $j=1,\dots,n$, one can define a path $\eta_j$ 
from $*$ to $-j+\frac{1}{4}$ constructed from a straight segment 
along the $x$-axis after replacing the segments joining 
$-k+\frac{1}{4}$ and $-k-\frac{1}{4}$ by $\delta_k^+$, $k=1,\dots,j-1$.

The following set of paths based on $*$ will be useful:
$x_j:=\eta_j\cdot\delta_j\cdot(\eta_j)^{-1}$, $j=1,\dots,n$.
It is well known that $x_1,\dots,x_n$ is a basis of the
free group $\pi_1(\bc\setminus\by^0; *)$. Moreover, 
$\gamma=x_n\cdot\ldots\cdot x_1$.
\medbreak

Let us consider the group $B_{\by^0}$. For $j=1,\dots,n-1$,
the braid $\sigma_j$ will be defined as a set 
$\{\gamma_1^{(j)},\dots,\gamma_n^{(j)}\}$
of paths such that:
\begin{itemize}
\item $\gamma_k^{(j)}$ is constant with value $-k$ if $k\neq j,j+1$;

\item $\gamma_j^{(j)}$ is $\delta_j^+$ and

\item $\gamma_{j+1}^{(j)}$ is $\delta_j^-$.
\end{itemize}

It is well known that $\sigma_1,\dots,\sigma_{n-1}$ generate
$B_{\by^0}$, providing a way to identify $B_{\by^0}$ and $B_n$.
These identifications are equivariant with respect
to the actions $\Phi$ and
$\Phi_{\by^0}$, which also become identified.
\end{ejm}
\bigbreak Based on this example we will describe 
isomorphisms from $\pi_1(\bc\setminus\by; *)$ to $\bff_n$ and 
from $B_{\by}$ to $B_n$ for any $\by\in X$, commuting with
respect to $\Phi$ and $\Phi_{\by}$. 
This can be achieved by defining braids with different ends.

\begin{dfn} Let $\by^1$ and $\by^2$ be subsets of $\bc$ with
exactly $n$ elements --\,that is, $\#\by^1=\#\by^2=n$. 
A \emph{braid from $\by^1$ to $\by^2$} is a homotopy class, 
relative to $\by^1 \cup \by^2$, of sets 
$\{\gamma_1,\dots,\gamma_n\}$ of paths $\gamma_j:[0,1]\to\bc$ 
starting at $\by^1$, ending at $\by^2$ and such that 
$\forall t\in[0,1]$, $\#\{\gamma_1(t),\dots,\gamma_n(t)\}=n$.
The sets $\by^1$ and $\by^2$ are usually referred to as the
\emph{ends of the braid}.
\end{dfn}

As above, we will usually not distinguish between a braid
and a representative of its class. The set of all braids from 
$\by^1$ to $\by^2$ is denoted by $B(\by_1,\by_2)$. 
Note that $B_{\by}=B(\by,\by)$ for any $\by\in X$.
Given $\by^1,\by^2,\by^3\in X$, there are natural definitions
\bdm
B(\by_1,\by_2)\times B(\by_2,\by_3)\to B(\by_1,\by_3),\quad
B(\by_1,\by_2)\to B(\by_2,\by_1),
\edm
for the product and inverse of braids, which have the natural
properties which identify them with the fundamental grupo\"{\i}d
of $X$ and will be called the 
\emph{braid grupo\"{\i}d on $n$ strings}.

Moreover, let us consider $\by^1,\by^2\in X$
and let $\Delta$ be a geometric disk
containing $\by^1\cup\by^2$ in its interior.

\begin{dfn} A \emph{$(\by^2,\by^1)$-special homeomorphism}
is an orientation-preserving homeomorphism $f:\bc\to\bc$ 
such that $f(\by^1)=\by^2$ and $f$ is the identity on 
$\bc\setminus\Delta$. Analogously, 
\emph{$(\by^2,\by^1)$-special isotopies} can be defined.
\end{dfn}

The quotient set of special isotopy classes of 
$(\by^2,\by^1)$-special homeomorphisms of $\bc$ will be
denoted by $\bbb(\by^2,\by^1)$. With the natural operations, 
one obtains a grupo\"{\i}d on $X$ which is naturally 
anti-isomorphic to the braid grupo\"{\i}d. By choosing
$*\in\partial\Delta$ as a base point, one has a natural 
mapping
\bdm
\bbb(\by^2,\by^1)\times\pi_1(\bc\setminus\by^1; *)\to
\pi_1(\bc\setminus\by^2; *)
\edm
which may be viewed as a grupo\"{\i}d action.
It induces a right grupo\"{\i}d action
\bdm
\pi_1(\bc\setminus\by^1; *)\times B(\by^1,\by^2)\to
\pi_1(\bc\setminus\by^2; *),
\edm
which becomes $\Phi_{\by}$ when $\by^1=\by^2=\by$.
This action will be denoted using exponential notation.

\medbreak
Let us fix ends $\by^1,\by^2\in X$ and a braid 
$\tau\in B(\by^1,\by^2)$. 
There are two natural isomorphisms
\bdm
\begin{array}{cccl}
R_\tau:&B_{\by^1}&\to&B_{\by^2}\\
&\sigma&\mapsto&\tau^{-1}\cdot\sigma\cdot\tau
\end{array}
\edm
and
\bdm
\begin{array}{cccl}
\Psi_\tau:&\pi_1(\bc\setminus\by^1; *)
&\to&\pi_1(\bc\setminus\by^2; *)\\
&x&\mapsto&x^\tau.
\end{array}
\edm
If $\by^0=\by^1=\by$, then $\tau\in B_\by$, 
$R_\tau$ is the inner automorphism 
$\sigma\mapsto\tau^{-1}\cdot\sigma\cdot\tau$
and $\Psi_\tau$ is the automorphism 
$\Phi_{\by}(\bullet,\tau)$ induced by the right
action $\Phi_{\by}$ with respect to~$\tau$.
We summarize these results.

\begin{prop} \label{otropunto} 
Let $\by^0=\{-1,\dots,-n\}$ as in 
\emph{example \ref{standard}}, and let $\by\in X$. 
We identify $\pi_1(\bc\setminus\by^0; *)$
with $\bff_n$ and $B_{\by^0}$ with $B_n$.

Then, for any $\tau\in B(\by^0,\by)$ there is a canonical
identification of $\pi_1(\bc\setminus\by; *)$
with $\bff_n$ and of $B_{\by}$ with $B_n$
by means of $\Psi_\tau$ and $R_\tau$ respectively.
These identifications are equivariant
with respect to $\Phi$ and $\Phi_{\by}$.
\smallbreak
Moreover, two such identifications differ
by an inner automorphism of $B_n$, with respect to 
a given $\hat\tau\in B_n$, and
the corresponding automorphism of $\bff_n$ is 
$\Phi(\bullet,\hat\tau)$.
\end{prop}

\section{Geometric bases and lexicographic braids}

The notations introduced along the previous section will 
also be used in the present one.
We recall an important definition from algebraic geometry.

\begin{dfn} Let $X$ be a connected projective manifold 
and let $H$ be a hypersurface of $X$. Let $*\in X\setminus H$ 
and let $K$ be an irreducible component of $H$. A homotopy 
class $\gamma\in\pi_1(X\setminus H;*)$ is called a 
\emph{meridian about $K$ with respect to $H$}
if $\gamma$ has a representative $\delta$ satisfying the
following properties: 
\begin{enumerate}[(a)]
\item there is a smooth disk 
$\Delta \subset X$ transverse to $H$ such that 
$\Delta\cap H=\{*'\} \subset K$ is a smooth point of $H$
\item there is a path $\alpha$ in $X\setminus H$ from $*$ 
to $*'' \in \partial\Delta$ 
\item $\delta=\alpha\cdot\beta\cdot\alpha^{-1}$, where
$\beta$ is the closed path obtained by traveling from $*''$ 
along $\partial\Delta$ in the positive direction.
\end{enumerate}
\end{dfn}

It is well known that any two meridians about $K$ with 
respect to $H$ are conjugate in $\pi_1(X\setminus H;*)$.
Moreover, the converse also holds.
\medbreak
In example \pref{standard}, $x_j$ is a meridian loop about $-j$
in $\pi_1(\bc\setminus\by^0;n+1)$
where $\by^0=\{-1,\dots,-n\}$, for $j=1,\dots,n$.
Let us note that if $\by^1\in X$ and $\tau\in B(\by^0,\by^1)$,
then $x_1^\tau,\dots,x_n^\tau$ are meridian loops about the 
points in $\by^1$, and of course they also form a basis for
the free group $\pi_1(\bc\setminus\by^1;*)$.

\begin{dfn}
\label{defgeom} 
Let us consider $\by$ a point of $X$, $\Delta\subset\bc$ a
geometric closed disk containing $\by$ in its interior,
and $*$ a point in $\partial\Delta$. Also, let 
$\gamma\in\pi_1(\bc\setminus\by,*)$ be as in 
\emph{notation \pref{vuelta}}. A \emph{geometric basis} of 
the free group $\pi_1(\bc\setminus\by,*)$
(\emph{with respect to $\Delta$})
is an ordered basis $\mu_1,\dots,\mu_n$ such that:
\begin{enumerate}
\item $\mu_1,\dots,\mu_n$ are meridians of
the points in $\by$;
\item $\mu_n\cdot\ldots\cdot\mu_1=\gamma$; in particular,
it is the inverse of a meridian about the point
at infinity of $\bc$.
\end{enumerate}
\end{dfn}

\begin{prop}
\label{orbgeom} 
Under the conditions and notations of 
\emph{example \pref{standard}} the following holds:
\begin{enumerate}[\rm(i)]
\item The base $x_1,\dots,x_n$
is geometric with respect to $\Delta$.
\item For any $\tau\in B_n$, the basis
$x_1^\tau,\dots,x_n^\tau$ is geometric with
respect to $\Delta$.
\label{propbas2}
\item Any geometric basis with respect to $\Delta$
is obtained as in \emph{\pref{propbas2}}.
\end{enumerate}
\end{prop}

Let us denote by $\cals(\Delta)$ the set of all geometric
bases with respect to $\Delta$.
The map $\Phi$ induces an action of $B_n$ on $\bff_n^n$
where $\cals(\Delta)$ becomes the orbit of $(x_1,\dots,x_n)$. 
The group $B_n$ acts freely on $\cals(\Delta)$.
The next result is a direct consequence of 
propositions \ref{otropunto} and \ref{orbgeom}.

\begin{cor} 
Let $\by\in X$, $\Delta$ and $*$ be as in
\emph{definition \pref{defgeom}}. Let us also assume that
$\Delta$ contains $\by^0$ in its interior. Then the set
$\cals(\Delta)$
is an orbit under the action of $B_{\by}$ induced
by $\Phi_{\by}$ on $\pi_1(\bc\setminus\by;*)^n$.

Note that, for any braid $\tau\in B(\by^0,\by)$, 
the set $x_1^\tau,\dots,x_n^\tau$ is a geometric
basis of $\pi_1(\bc\setminus\by;*)$. Moreover, any
element of the orbit can be obtained in this manner.
\end{cor}

Our next purpose is to construct canonical 
isomorphisms between $B_{\by}$ and $B_n$
for any $\by\in X$. This can be done by 
fixing a special braid $\tau\in B(\by^0,\by)$.
We must emphasize that even if the process
to produce such a braid is canonical, this
method is not invariant by orientation-preserving
homeomorphisms of $\bc$ (not even by rotations).
\medbreak

\begin{cnt}
\label{lexico-const}
We will start by constructing lexicographic
bases and braids. Let us fix $\by\in X$ and let
us consider $\by^0$ as above. Let $\Delta$ be
a geometric disk (centered at $0$) containing 
$\by\cup\by^0$ and choose $*$ as the only positive 
real number in $\partial\Delta$.
\smallbreak
Let us consider the basis $x_1,\dots,x_n$ of 
$\pi_1(\bc\setminus\by^0,*)$ given in example 
\ref{standard}.
\smallbreak
One can order the points $t_1,\dots,t_n\in\by$ 
such that $\Re t_1,\geq\dots\geq\Re t_n$ and if 
$\Re t_j=\Re t_{j+1}$, then $\Im t_j>\Im t_{j+1}$.
This ordering is called the 
\emph{lexicographic ordering} of $\by$.

Choose $\varepsilon>0$ such that the closed disks 
centered at $t_j$, $j=1,\dots,n$ of radius 
$\varepsilon$ are pairwise disjoint and
contained in the interior of $\Delta$.
Let us consider the polygonal path $\Gamma$ joining
$*,t_1,\dots,t_n$ and the circles of radius 
$\varepsilon$ centered at $t_1,\dots,t_n$. For each
$j=1,\dots,n$, we construct $y_j$, a meridian about 
$t_j$ with respect to $\bc$ based at $*$. For this 
purpose, we will follow the process shown in example 
\pref{standard}, but replacing the segment on the real 
line by $\Gamma$ and the disks of radius $\frac{1}{4}$ 
centered at $-j$ by the disks of radius $\varepsilon$ 
centered at $t_j$, $j=1,\dots,n$.
\smallbreak

In this way a geometric basis $y_1,\dots,y_n$
of $\pi_1(\bc\setminus\by;*)$ is produced. Such
a basis will be called the \emph{lexicographic basis} 
for $\by$. Note that, in particular, $x_1,\dots,x_n$ 
is the lexicographic basis for $\by^0$.

Note that there is a unique
braid $\tau_\by\in B(\by_0,\by)$ such that
$y_j=x_j^{\tau_\by}$, $\forall j=1,\dots,n$. This 
braid $\tau_\by$ is called the 
\emph{lexicographic braid} associated to $\by$.
The corresponding isomorphism 
$R_\by:=R_{\tau_\by}:B_n\to B_\by$
will be called \emph{lexicographic isomorphism}
from $B_n$ to $B_\by$.
\end{cnt}

\begin{obs}
\label{lexico} 
Let us consider $\by_1,\by_2\in X$.
Lexicographic braids allow us to consider
canonical bijections of $B(\by_1,\by_2)$
where $B_n=B_{\by^0}$:
\bdm
\begin{array}{ccl}
B(\by^1,\by^2)&\to&B_n\\
\sigma&\mapsto&\tau_{\by^1}\cdot\sigma\cdot\tau_{\by^2}^{-1}.
\end{array}
\edm
Let us say we have a braid in $\bc\times[0,1]$, then a generic 
projection onto $\br\times[0,1]$ is obtained as follows:
\begin{itemize}

\item 
Take the projection $\Re:\bc\to\br$ given by the real part.

\item 
For any isolated double point of the projection, 
the branch with the smallest imaginary part will be drawn 
using a continuous line.

\item 
If a line of double points occurs (e.g, a pair of conjugate
arcs), we  slightly perturb the projection in order
to have greater imaginary parts to the right and smaller 
imaginary parts to the left.
\end{itemize}
\end{obs}

\section{Braid monodromy and fibered curves}

In order to define braid monodromy and fibered curves,
we will consider the affine plane $\bc^2$ and the
projection $\pi:\bc^2\to\bc$ given by $\pi(x,y)=x$.
Instead of fixing the second coordinate of $\bc^2$ 
we will allow affine transformations, preserving $\pi$,
to act on $\bc^2$.

\begin{dfn} 
A reduced affine curve $\cc\subset\bc^2$ is said to be
\emph{horizontal} with respect to $\pi$ if 
$\pi_{|\cc}:\cc\to\bc$ is a proper map. 
The horizontal degree $\deg_\pi(\cc)$ of $\cc$
is the degree of $\pi_{|\cc}$, i.e., the number of 
preimages of a regular value. The set of regular values 
of $\pi_{|\cc}$ will be denoted by $\bc_\cc$.
\end{dfn}

Let $\cc$ be a horizontal curve of horizontal degree 
$d$. Note that, after fixing the second coordinate of 
$\bc^2$ we can assume that an equation $f(x,y)=0$ for
$\cc$ is given by
\begin{displaymath}
f(x,y)=y^d+f_1(x) y^{d-1}+\dots+f_1(x)y+f_0(x),
\end{displaymath}
where $f_j(x)\in\bc[x]$, $j=0,1,\dots,d-1$.
The condition of being horizontal is equivalent
to the non-existence of vertical asymptotes
(vertical lines included).
Note that the set of critical values is the set of 
zeros of the discriminant of $f$ with respect to $y$, 
which is a polynomial $\fd(x)\in\bc[x]$.
Denoting this set of zeros by $x_1,\dots,x_r$, one 
has $\bc_\cc=\bc\setminus\{x_1,\dots,x_r\}$. Therefore 
the special fibers of $\pi$ are exactly the vertical 
lines $L_i:=\pi^{-1}(x_i)$ of equation 
$x=x_i$, $i=1,\dots,r$.

\begin{dfn} Let $\cc\subset\bc^2$ be a
horizontal curve with respect to $\pi$.
With the above notation, the \emph{fibered curve}
associated to $\cc$ is the curve
$\cc^\varphi:=\cc\cup L_1\cup\dots\cup L_r$.
\end{dfn}

The motivation behind this definition is the following.
Let $\check \cc:=\pi_{|\cc}^{-1}(\bc_\cc)$.
Since $\pi_{|\cc}$ is proper, 
$\tilde\pi:\check \cc\to\bc_\cc$ is a 
(possibly non-connected) covering map. The mapping
$\pi^\varphi:\bc^2\setminus \cc^\varphi\to\bc_\cc$ 
is a locally trivial fibration whose fiber is 
diffeomorphic to $\bc\setminus\{1,\dots,d\}$.
The polynomial $f$ induces an algebraic mapping
$\tilde f:\bc_\cc\to X$ defined as 
$\tilde f(x_0):= \{x=x_0\} \cap \cc$ 
--\,or equivalently as the polynomial 
$f(x_0,t) \in \bc[t]$.
Let us fix $*\in\bc_\cc$, a complex regular value 
on the boundary of a geometric disk $\Delta$
containing $\{x_1,\dots,x_r\}$ in its interior. 
Let us denote by $\by^*$ the set of roots
of the polynomial $f(*,t)$.

\begin{dfn} The homomorphism 
$\nabla_*:\pi_1(\bc_\cc;*)\to B_{\by^*}$ induced 
by $\tilde f$ is called a \emph{braid monodromy} 
of $\cc$ with respect to $\pi$.
\end{dfn}

\begin{obs} 
Note that $\nabla_*$ classifies the locally 
trivial fiber bundle $\pi^\varphi$.
\end{obs}

Let us briefly explain how to construct this braid 
monodromy. Note that there is a particular class of 
basis for $\pi_1(\bc_\cc;*)$, namely the geometric 
bases with respect to $\Delta$. Let us fix one of 
these bases, say $\mu_1,\dots,\mu_r$ 
(e.g., one can choose the lexicographic basis). 
Let us consider the lexicographic
isomorphism $R_{\by^*}:B_n\to B_{\by^*}$.
These facts allow us to represent $\nabla_*$ by
an element $(\tau_1,\dots,\tau_r)\in B_d^r$ such
that:
\bdm
\tau_j:=
R_{\by^*}^{-1}(\nabla_*(\mu_j)),\quad \forall j=1,\dots,r.
\edm

\begin{dfn} 
We say that $(\tau_1,\dots,\tau_r)\in B_d^r$
\emph{represents the braid monodromy of $\cc$} if
there exist:
\begin{itemize} 

\item a geometric disk $\Delta$ containing
$\{x_1,\dots,x_r\}$ in its interior;

\item an element $*\in\partial\Delta$ such that
$\by^*$ is the set of roots of the polynomial 
$\tilde f(*) \in \bc[t]$;

\item a geometric basis $\mu_1,\dots,\mu_r$ of
$\pi_1(\bc_\cc,*)$;

\item a braid $\beta\in B(\by^*,\by^0)$
\end{itemize}
for which $\tau_j=R_{\beta}(\nabla_*(\mu_j))$,
$\forall j=1,\dots,r$.
\end{dfn}

\begin{obs} 
Two natural right actions on $B_d^r$ are related 
to the concept of $r$-tuples of braids representing 
the monodromy.

The first one is an action produced by $B_r$ as 
follows: 
\begin{displaymath}
(\tau_1,\dots,\tau_r)^{\sigma_i}:=
(\tau_1,\dots,\tau_{i-1},\tau_{i+1},\tau_{i+1}
\tau_i\tau_{i+1}^{-1},\tau_{i+2},\dots,\tau_r),
\end{displaymath}
where $\sigma_i\in B_r$, $i=1,\dots,r-1$ 
is a canonical generator
and $(\tau_1,\dots,\tau_r)\in B_d^r$.
The second action is given by $B_n$ and it is 
defined by conjugation on each coordinate as follows:
\begin{displaymath}
(\tau_1,\dots,\tau_r)^{\beta}:=
(\tau_1^{\beta},\dots,\tau_r^{\beta}),
\end{displaymath}
where $\beta\in B_d$ and $(\tau_1,\dots,\tau_r)\in B_d^r$.
These actions commute with each other and hence they 
define a new action of $B_r\times B_d$.
The first action represents the change of
geometric basis. The second action represents both
the change of the chosen braid in $B(\by^*,\by^0)$ and 
also the change of the base point.
\end{obs}

This can be summarized as follows:

\begin{prop}
\label{rep} 
Let $(\tau_1,\dots,\tau_r)$ be an $r$-tuple of braids 
representing the braid monodromy of $\cc$ w.r.t. $\pi$. 
Then $(\tilde\tau_1,\dots,\tilde\tau_r)\in B_d^r$
represents the braid monodromy of $\cc$ if and only 
if $(\tilde\tau_1,\dots,\tilde\tau_r)$
is in the orbit of $(\tau_1,\dots,\tau_r)$ by the
action of $B_r \times B_d$ on $B_d^r$.
\end{prop}

\begin{obs} 
Moreover, note that cyclic permutations of
$r$-tuples representing a braid monodromy
have essentially the same information as their
original representatives. This remark leads to
the concept of \emph{pseudogeometric bases}, 
that is, a basis $\mu_1,\dots,\mu_r$ such that 
$\mu_j$ is a meridian about $x_{\sigma(j)}$ for 
some $\sigma\in\Sigma_r$ and 
$(\mu_r\cdot\ldots\cdot\mu_1)^{-1}$ is a
meridian about the point at infinity.
Sometimes, braid monodromies are easier to 
compute for pseudogeometric bases.
\end{obs}

\begin{dfn} 
Two affine horizontal curves are said to have 
\emph{equivalent braid monodromies} if they have 
the same representatives for their braid monodromies.
\end{dfn}

Note that this equivalence relation is finer than 
the one arising from curves with the \emph{same} 
monodromy representation.

In general, it is difficult to find effective 
invariants to compare braid monodromies. For example, 
the invariants suggested by Libgober depend 
essentially on the conjugation class of the image 
of the braid monodromy. We recall how braid 
monodromy is related to the fundamental group.

\begin{zvk1}\label{zvk1} 
Let $\bff_d$ be the free group generated by a 
geometric basis $g_1,\dots,g_d$ and let 
$(\tau_1,\dots,\tau_r)\in B_d^r$ be an $r$-tuple 
of braids representing the braid monodromy of $\cc$.
Then the fundamental group of $\bc^2\setminus \cc$
is isomorphic to
\begin{displaymath}
\left\langle
g_1,\dots,g_d\big |
g_j^{\tau_i}=g_j,\quad
i=1,\dots,r,\quad
j=1,\dots,d
\right\rangle,
\end{displaymath}
\end{zvk1}

The main tools required to prove this theorem are
the classical Van Kampen theorem and this fibered
version of Zariski-Van Kampen's theorem.

\begin{prop}\label{zvk2}
The fundamental group of $\bc^2\setminus \cc^\varphi$
is isomorphic to
\begin{displaymath}
\left\langle
g_1,\dots,g_d,\alpha_1,\dots,\alpha_r\big |
g_j^{\tau_i}=\alpha_i^{-1}g_j\alpha_i,\quad
i=1,\dots,r,\quad
j=1,\dots,d
\right\rangle.
\end{displaymath}
\end{prop}

The key point in the proof of this fibered version 
is the long exact sequence of homotopy associated 
with a fibration. 

More information about the pair $(\bc^2,\cc^\varphi)$ 
can be obtained by choosing the paths representing 
$g_1,\dots,g_d,\alpha_1,\dots,\alpha_r$ in an 
appropriate manner.

\begin{dfn} 
Let $\mu_1,\dots,\mu_r$ be a geometric basis of
$\pi_1(\bc_\cc,*)$. We assume that these paths 
have their support in a simply connected compact 
set $K\subset\bc$ satisfying $*\in\partial K$.
Let $L\subset\bc^2$ be a simply connected compact 
subset of $\bc^2$ such that 
$\pi^{-1}(\cc\cap K)\subset L$ and let
$\hat *\in\bc$ such that $(*,\hat *)\notin L$. 
We say that 
$\alpha_1,\dots,\alpha_r\in
\pi_1(\bc^2\setminus \cc^\varphi;(*,\hat *))$
is a \emph{suitable lifting} of $\mu_1,\dots,\mu_r$ if
\begin{enumerate}[(1)]

\item $\pi^\varphi_*(\alpha_j)=\mu_j$, $j=1,\dots,r$;

\item the support of $\alpha_1,\dots,\alpha_r$ is
in $\pi^{-1}(K)\setminus L$;

\item for any $j=1,\dots,r$, the closed path 
$\alpha_j$ is a meridian about the line $L_{\sigma(j)}$ 
with respect to $\bc^2$, where $\sigma$ is the 
permutation associated to the original geometric basis.
\end{enumerate}
\end{dfn}

\begin{lema} 
Let $\mu_1,\dots,\mu_r$ be a geometric basis of
$\pi_1(\bc_\cc,*)$. Then, a suitable lifting of this
basis is unique up to homotopy.
\end{lema}

\begin{prop}\label{zvk3}
Let $L^\varphi_*:=(\pi^{\varphi})^{-1}(*)$ be a generic 
fiber of $\pi^{\varphi}$ and 
let $g_1,\dots,g_d$ be a geometric basis of 
$\pi_1(L^\varphi_*,(*,\hat *))$. Then both the elements 
$g_1,\dots,g_d$ and a suitable lifting 
$\alpha_1,\dots,\alpha_r$ of $\mu_1,\dots,\mu_r$ may be 
chosen for the presentation of
$\pi_1(\bc^2\setminus \cc^\varphi;(*,\hat *))$ 
in theorem \ref{zvk2}.
\end{prop}

\begin{cor} Any presentation of 
$\pi_1(\bc^2\setminus \cc^\varphi;(*,\hat *))$ in terms
of a geometric basis of $\pi_1(L^\varphi_*,(*,\hat *))$ and
a suitable lifting $\alpha_1,\dots,\alpha_r$ of
a geometric basis of $\pi_1(\bc_\cc;*)$, determines the 
braid monodromy $\Phi$ of $\cc$.
\end{cor}

The following definition and lemma will help to give 
a more geometrical construction for the suitable 
lifting described in proposition \ref{zvk3}. 

\begin{dfn} 
Let $A,B\gg 0$. A polydisk 
$\tilde\Delta:=\Delta_A\times\Delta_B\subset\bc^2$ of 
multiradius $(A,B)$ is \emph{well adapted} to $\cc$ if 
$\{x_1,\dots,x_r\}$ is contained in the interior of 
$\Delta_A$ and
\bdm
\{(x,y)\in \cc|\  x\in\Delta_{A}\}
\subset\Delta_{A}\times\Delta_{B}.
\edm
\end{dfn}

\begin{lema} 
There exists a real number $A_0>0$ such that for any 
$A\geq A_0$ there exists $B_0(A)>0$ satisfying that 
for any $B\geq B_0(A)$ the polydisk 
$\tilde\Delta:=\Delta_A\times\Delta_B$ 
is well adapted to $\cc$.
\end{lema}

\begin{obs} If $\tilde\Delta$ is a well-adapted polydisk
as above, we can set $(*,\hat*)=(A,B)$ and choose the
paths $\alpha_i$ in $\Delta_A\times\{B\}$.
\end{obs}

Now we can state the main theorem relating equivalence
of pairs of fibered curves and equivalence of braid
monodromies. This is a partial converse to \cite{kt:00}
and \cite{car:xx}.

\begin{thm}
\label{top} 
Let $\cc_1,\cc_2\subset\bc_2$ be two horizontal affine
curves and let us consider the standard embedding
$\bc^2\subset\bp^2$. If
$F:(\bc^2,\cc_1^\varphi)\to(\bc^2,\cc_2^\varphi)$,
is an orientation preserving homeomorphism
that extends to a homeomorphism of $\bp^2$,
then $\cc_1$ and $\cc_2$ have equivalent braid monodromies.
\end{thm}

\begin{proof}
For the proof we will use the notation introduced in
statements \ref{zvk1}, \ref{zvk2} and \ref{zvk3}.

Let us consider a horizontal curve $\cc$.
The elements $g_1,\dots,g_d$ freely generate 
a subgroup $H$, which is normal since 
$H=\ker\pi^\varphi_*$.
Note that $g_1,\dots,g_d$ are meridians about (all) 
the (irreducible) components of $\cc$.
Since $H$ is normal, it is generated by all the 
meridians about all irreducible components of $\cc$.
The group $H$ is, hence, identified with the fundamental 
group of the generic fiber. Therefore, the basis 
$g_1,\dots,g_d$ is geometric. Note that 
$g_d\cdot\ldots\cdot g_1$ is the boundary of a 
sufficiently big disk on the fiber.
One obtains a natural exact sequence:
\begin{equation}
\label{sucex}
1\to H\longrightarrow
\pi_1(\bc^2\setminus \cc^\varphi;(*,\hat *))
\overset{\pi^\varphi_*}{\longrightarrow}
\pi_1(\bc\setminus \bc_\cc;*)\to 1.
\end{equation}
\smallbreak

Given any two curves $\cc_1,\cc_2$ as in the statement,
let us take the elements $g_1^{(1)},\dots,g_d^{(1)}$ and
$\alpha_1^{(1)},\dots,\alpha_r^{(1)}$ described in 
proposition \ref{zvk3}. We recall that the elements
$g_1^{(1)},\dots,g_d^{(1)}$ form a geometric basis
of $H^{(1)}$ and ${\pi^\varphi_*}^{(1)}(\alpha_1^{(1)}),
\dots,{\pi^\varphi_*}^{(1)}(\alpha_r^{(1)})$
is a geometric basis of $\pi_1(\bc\setminus \bc_{\cc_1};*)$.
Let us choose a big disk $\Delta_{A_1}$ of radius $A_1\gg 0$,
such that the values of the non-transversal vertical
lines to $\cc_1$ are in the interior of $\Delta_{A_1}$.
Note that we can construct a polydisk $\tilde\Delta_1$
well adapted to $\cc_1$ and two polydisks $\tilde\Delta_2$
and $\hat\Delta_2$ well adapted to $\cc_2$ such that 
\bdm
\tilde\Delta_2\subset F(\tilde\Delta_1)\subset\hat\Delta_2.
\edm
Let us define
\bdm
g_j^{(2)}:=F_*(g_j^{(1)}),\quad j=1,\dots,d.
\edm
From the discussion above, $g_1^{(2)},\dots,g_d^{(2)}$ is a
basis for the free group $H^{(2)}$. Moreover, let us suppose
that the generic fiber $L^{(1)}_*$ (of $\pi$) is very close 
to a non-transversal vertical line $L^{(1)}_i$ (if no such 
line exists there is nothing to prove). Since the boundary 
of a big disk in $L^{(1)}_i$ is sent to the boundary
of a big disk in $L^{(2)}_i$, we easily deduce that 
$g_d^{(2)}\cdot\ldots\cdot g_1^{(2)}$ is homotopic to the 
boundary of a big disk in a generic fiber for
$\pi^{\varphi(2)}$. This argument proves that 
$g_1^{(2)},\dots,g_d^{(2)}$ is a geometric basis of $H^{(2)}$.

\smallbreak
Let us define
\bdm
\alpha_i^{(2)}:=F_*(\alpha_i^{(1)}), \quad i=1,\dots,r.
\edm
By the naturality of the exact sequence \pref{sucex}
we deduce that ${\pi^\varphi_*}^{(2)}(\alpha_1^{(2)})$, \dots,
${\pi^\varphi_*}^{(2)}(\alpha_r^{(2)})$
is a basis of the free group $\pi_1(\bc\setminus \bc_{\cc_2};*)$.
Since the extension of $F$ to $\bp^2$ preserves the line at 
infinity and since the inverse of 
${\pi^\varphi_*}^{(1)}(\alpha_r^{(1)}
\cdot \ldots\cdot\alpha_1^{(1)})$ is a meridian about the point
at infinity in $\bc$, we deduce that it is also the case
for ${\pi^\varphi_*}^{(2)}(\alpha_r^{(2)}\cdot
\ldots\cdot\alpha_1^{(2)})^{-1}$.
Then, ${\pi^\varphi_*}^{(2)}(\alpha_1^{(2)}),
\dots,{\pi^\varphi_*}^{(2)}(\alpha_r^{(2)})$ is a 
pseudogeometric basis of $\pi_1(\bc\setminus \bc_{\cc_2};*)$.
\smallbreak
From the exact sequence \pref{sucex} and the system
of generators $g_i^{(2)},\alpha_j^{(2)}$, we obtain a 
representation of $\bff_r$ in $B_d$ which provides an
element of $M_2\in (B_d)^r$. We have proven that this
element represents the braid monodromy of $\cc_2$.

\smallbreak
Because of the way the generators are chosen, one can apply the 
same argument to $\cc_1$, obtaining an element $M_1\in (B_d)^r$ 
representing the braid monodromy of $\cc_1$. Since the two 
families of generators are related by $F_*$, we deduce 
that $M_1=M_2$
\end{proof}

This theorem will be used to compare curves in $\mm_1$
and $\mm_2$. We must compute their braid monodromies
and find suitable invariants in order to
be able to compare them.

\section{Effective invariants of braid monodromy}
\label{efinv}

Let $G$ be a group and let $r$ be a positive integer.
Let us consider the sets $G^r$ and $G\backslash G^r$
which is the quotient of $G^r$ by the diagonal action of
$G$ by conjugation. The braid group $B_r$ acts on $G^r$
by the so-called Hurwitz action:
\bdm
\sigma_j\cdot(g_1,\dots,g_r):=(g_1,\dots,g_{j-1},g_{j+1},
g_{j+1}g_jg_{j+1}^{-1},g_{j+2},\dots,g_r),\ 1\leq j<r,
\edm
where $g_1,\dots,g_r\in G$. Since Hurwitz and 
conjugation actions commute, $B_r$ also acts on 
$G\backslash G^r$.
\smallbreak

Let us denote by $\nn_r(G)$ (resp. $\mm_r(G)$)
the quotient of $G^r$ (resp. $G\backslash G^r$)
by the Hurwitz action. The index $r$ will be dropped
if no ambiguity seems likely to arise. The elements of
$\mm_r(G)$ will be called $G$-monodromies (of order $r$).
Before we use the main result to produce effective 
invariants, let us set the main definitions of this last 
section.

\begin{dfn}
\label{bmb} 
Let $\cc$ be a horizontal curve with $\deg_\pi(\cc)=d$ 
and possessing $r$ non-transversal vertical lines. 
The \emph{braid monodromy} of $\cc$ is defined as the element 
in $\mm_r(B_d)$ determined by any $r$-tuple representing the
braid monodromy of $\cc$.
\end{dfn}

\begin{obs} 
This point of view is inspired by the work of Brieskorn
on \emph{automorphic sets} \cite{brk:88}. 
We restrict our attention to automorphic sets defined by 
conjugation on groups. We have also modified the conjugation 
action defined in \cite{brk:88}.
\end{obs}

\begin{dfn} 
Let $\gog:=(g_1,\dots,g_r)\in G^r$. The 
\emph{pseudo Coxeter element} associated with $\gog$ is 
defined as $c(\gog):=g_r\cdot\ldots\cdot g_1$.
\end{dfn}

Note that pseudo Coxeter element is also well defined in 
$\nn(G)$ and its conjugation class is well defined in 
$\mm(G)$.

Let $\phi:G_1\to G_2$ be a group homomorphism. It induces
in a functorial way mappings
$\phi_{\nn}:\nn(G_1)\to\nn(G_2)$ and
$\phi_{\mm}:\mm(G_1)\to\mm(G_2)$.

\begin{dfn}
\label{gmb} 
Let $\cc$ be as in definition \ref{bmb}
and let $\Phi:B_d\to G$ be a representation of $B_d$ onto
a group $G$. Then the \emph{$(G,\Phi)$-monodromy} of $\cc$ 
is the image in $\mm_r(G)$ by $\Phi_{\mm}$ of the braid 
monodromy of $\cc$.
\end{dfn}

\begin{prop} 
Let $\cc_1,\cc_2$ be horizontal curves with 
$\deg_\pi(\cc_1)=\deg_\pi(\cc_2)=d$ having
$r$ non-transversal vertical lines. Let $\Phi:B_d\to G$ 
be a representation. If the $(\Phi,G)$-monodromies of 
$\cc_1$ and $\cc_2$ are not equal then there is no 
orientation-preserving homeomorphism of $(\bc^2,\cc_1)$
and $(\bc^2,\cc_2)$ that extends to a homeomorphism of 
$\bp^2$.
\end{prop}

\begin{proof} 
This is a straightforward consequence of theorem \ref{top}.
\end{proof}

If $G$ is a finite group then $\mm(G)$ is a
finite set and hence, knowing the braid monodromies of
two given curves would allow us to compare their 
$(\Phi,G)$-monodromies, up to computational capacity. 
In the Appendix we give an algorithm implemented on GAP4 
\cite{GAP4}. Its input consists of braid monodromies of 
two curves and a finite representation of the braid group.
Its output affirms or negates the equality of their 
$(\Phi,G)$-monodromies.
We sketch the general lines of the algorithm:

\begin{enumerate}[{\bf (i)}]
\smallbreak\item Compute the $(\Phi,G)$-monodromies of 
the curves $\cc_1,\cc_2$.
\smallbreak\item Compute their pseudo Coxeter elements 
$c_G(\cc_1),c_G(\cc_2)\in G$. If they are not conjugate, 
braid monodromies are not equal. If they are conjugate
to each other, choose an element $g\in G$ such that 
$c_G(\cc_2)^g=c_G(\cc_1)=:h$. Let $H$ be the centralizer 
of $h$ in $G$ and consider the set $H_{\cc_2}$ of 
conjugates of the $(\Phi,G)$-monodromy
of $\cc_2$ by $gx$, $x\in H$. 
\smallbreak\item Compute the orbit of the 
$(\Phi,G)$-monodromy of $\cc_1$ by the action of $B_r$. 
Note that since Hurwitz action preserves pseudo Coxeter 
elements it is enough to consider conjugation by $H$. 
Since $B_r$ admits a generator system with two 
elements, an algorithm can be easily programmed to 
construct the finite orbit. For each new element of 
the orbit, verify if it is in $H_{\cc_2}$ and in that
case stop the program.
\item If no element of the orbit is in $H_{\cc_2}$ then
the two $(\Phi,G)$-monodromies are not equivalent.
\end{enumerate}

\begin{str}
\label{pasos} 
In general, in order to distinguish braid monodromies
for two horizontal curves $\cc_1$ and $\cc_2$ with
$\deg_\pi(\cc_1)=\deg_\pi(\cc_2)=d$, 
we proceed as follows:
\begin{itemize}

\item Verify if the curves (with the projections) 
have the same combinatorics (also at infinity).

\item If this is the case, then compute the fundamental
group of the curves.

\item If either the groups are isomorphic or we cannot
determine that they are not, then we compute the image
of the braid monodromy.

\item If we cannot determine whether or not the images are 
conjugate, we look for Libgober invariants which provide
easy-to-compare polynomials. Also the sequence of 
characteristic varieties might help to distinguish the
groups.

\item If the previous stpes do not work, 
we try the methods described in this section.
\end{itemize}
\end{str}

\section{Finite representations of the braid group}
\label{finitrep}

Finding finite representations of braid groups
is an interesting problem already studied by
several authors: Assion \cite{ass:78}, Kluitmann
\cite{kl:88}, Birman-Wajnryb \cite{brwj:86}, Wajnryb
\cite{waj:88,waj:91}. 
Infinite families of presentations have been obtained,
for instance the isomorphism $B_3\to SL(2;\bz)$ 
produces finite representations on $SL(2;\bz/n \bz)$.
These presentations can be carried over $B_4$ via the 
epimorphism $B_4\to B_3$. Analogously, homomorphisms 
onto symplectic groups provide finite representations.

\smallbreak
As it is well known, for a group $G$,
Hurwitz actions on $G^r$ and $G\backslash G^r$
provide finite representations of $B_r$.
Let $\gog\in G^r$ and let us denote by
$\Omega_\gog$ its orbit. Then the
Hurwitz action defines a homomorphism
$\tilde\theta_\gog:B_r\to\Sigma_{\Omega_\gog}$. 
Let us denote by $G_\gog$ the image 
$\tilde\theta_\gog(B_d)$, which is a transitive 
subgroup of $\Sigma_{\Omega_\gog}$. The induced 
mapping $\theta_\gog:B_d\to G_\gog$ is a surjective 
representation of $\gog$.

\smallbreak
Analogously, considering the class $[\gog]$ of 
$\gog$ in $G\backslash G^r$, one can construct another 
finite representation $\theta_{[\gog]}:B_d\to G_{[\gog]}$ 
which factors through the canonical mapping 
$G_\gog\to G_{[\gog]}$.

\begin{ejm} 
Let us suppose that $g_1,\dots,g_d$ commute pairwise.
Then, Hurwitz action factors through the canonical
map $B_d\to\Sigma_d:=\Sigma_{\{1,\dots,d\}}$ onto the
permutation action of $\Sigma_d$ on the ordered coordinates 
of $\gog$. If $g_1,\dots,g_d$ are pairwise distinct, 
then $G_\gog$ is naturally isomorphic to $\Sigma_d$.
\end{ejm}

In the general case, one can understand this braid 
action as a lifting of the permutation action on the 
abelianized group of $G$. In order to see an easy, but 
not trivial, example we consider the case $r=3$ and $G=\Sigma_3$.

Let us take $\gog:=[(1,2),(1,3),()]$.
Using GAP, we find that $\#\Omega_\gog=9$ and
$G_\gog$ is a group of order $162$ given by the
following presentation:
\bdm
\left\langle a,b\ :\
aba=bab,\ a^6=1,\ [a^2,b^2]=1,\ (ab^{-1})^3=1 
\right\rangle.
\edm
Note that only three relations (including the first 
two) are needed.

For $\gog:=[(1,2),(1,3),(1,2)]$, we find that 
$\#\Omega_\gog=8$ and $G_\gog$ is a group of order $24$ 
given by the following presentation:
\bdm
\left\langle a,b\ :\
aba=bab,\ a^3=1
\right\rangle.
\edm

For $\gog:=[(1,2),(1,2,3),()]$ and 
$\gog:=[(1,2,3),(1,2,3),(1,2)]$, 
we find that $\#\Omega_\gog=12$ and
$G_\gog$ is a group of order $48$ given by the 
following presentation:
\bdm
\left\langle a,b\ :\
aba=bab,\ a^4=1,\ (a^2b^{-2})^2=1,\ (ab^{-1})^3=1
\right\rangle.
\edm
As in the first case, only three relations including 
the first two are needed.

For $\gog:=[(1,2),(1,2),(1,2,3)]$, we find that 
$\#\Omega_\gog=18$ and $G_\gog$ is a group of order 
$17496$. Its presentation is more complicated in this case.

For $\gog:=[(1,2),(1,3),(1,2,3)]$, we find that 
$\#\Omega_\gog=18$ and $G_\gog$ has order $648$.

In this way, we have obtained all the $G_\gog$,
$\gog\in G^3$, up to conjugation.

\section{Construction of curves in $\mm$}

We follow the ideas in \cite{acct:00} to construct curves 
in $\mm$. Any curve in $\mm$ is projectively equivalent 
to exactly one of the following two projective curves 
$\tilde \cc_\beta:=\{f_{\beta}g_{\beta}=0\}$ with:
\begin{multline*}
f_{\beta}(x,y,z):=
{x}^{2}{y}^{3}+\left (303-216\,\beta\right )x{y}^{2}{z}^{2}+
\left (-636+450\,\beta\right )xy{z}^{3}+\\
+\left (-234\,\beta+331\right )x{z}^{4}
+\left (-18\,\beta+27\right )y{z}^{4}+\left (18\,\beta-26
\right ){z}^{5}\\
g_{\beta}(x,y,z):=x+\left ({\frac {10449}{196}}-
{\frac {3645}{98}}\,\beta\right )y+\left (-
{\frac {432}{7}}+{\frac {297}{7}}\,\beta\right )z
\end{multline*}
where $\beta^2=2$. The line $y=0$ is the tangent line to both 
curves $\tilde \cc_\beta$ at the $\be_6$ point which is $[1:0:0]$. 
We take the affine plane of coordinates $(x,z)$ and consider 
the projection $\pi(x,z)=z$. Since 
${\rm mult}_{\be_6}(\cc_\beta)=3$, one has 
$\deg_\pi(\cc_\beta)=3$. There are $4$ non-transversal vertical
lines corresponding to the singular points of types
$\ba_7$, $\ba_3$, $\ba_2$, $\ba_1$. There is also an ordinary 
tangent vertical line intersecting $\tilde \cc_\beta$ at a point 
of tangency denoted by $\ba_0$. The values for $z$ at these lines
are shown in table \ref{tablavert}
\begin{table}[ht]
\begin{center}
\begin{tabular}{|c|c|c|c|c|}\hline
$\ba_7$ & $\ba_3$ & $\ba_2$& $\ba_1$&$\ba_0$\\
\hline
$\frac{90+9\beta}{98}$ & $0$ & $1$ & $\frac{18+27\beta}{56}$ &%
 $-\frac{45+36\beta}{7}$ \\ \hline
\end{tabular}
\end{center}
\caption{}
\label{tablavert}
\end{table}

Let us denote by $\cc_{\beta}$ the affine curve 
$f_{\beta}(x,1,z)g_{\beta}(x,1,z)=0$. In both
cases we obtain $\bc_{\cc_{\beta}}$ by eliminating
five real points of $\bc$. We choose $*\in\br$, $*\gg 0$ and we 
choose geometric bases for $\pi_1(\bc_{\cc_\beta};*)$ using
the lexicographic construction \ref{lexico-const} 
and remark \ref{lexico}, with respect to a
segment in the real axis. We will denote by $\alpha_j^\beta$
the meridian about the non-generic line passing through the 
point $\ba_j$. By computing numerical values we obtain that the
geometric bases are
\bdm
\alpha_7^{\sqrt{2}},\alpha_1^{\sqrt{2}},
\alpha_2^{\sqrt{2}},\alpha_3^{\sqrt{2}},
\alpha_0^{\sqrt{2}}
\edm
and
\bdm
\alpha_2^{-\sqrt{2}},\alpha_0^{-\sqrt{2}},
\alpha_7^{-\sqrt{2}},\alpha_3^{-\sqrt{2}},
\alpha_1^{-\sqrt{2}}.
\edm

\begin{figure}[ht]
\includegraphics[height=3cm]{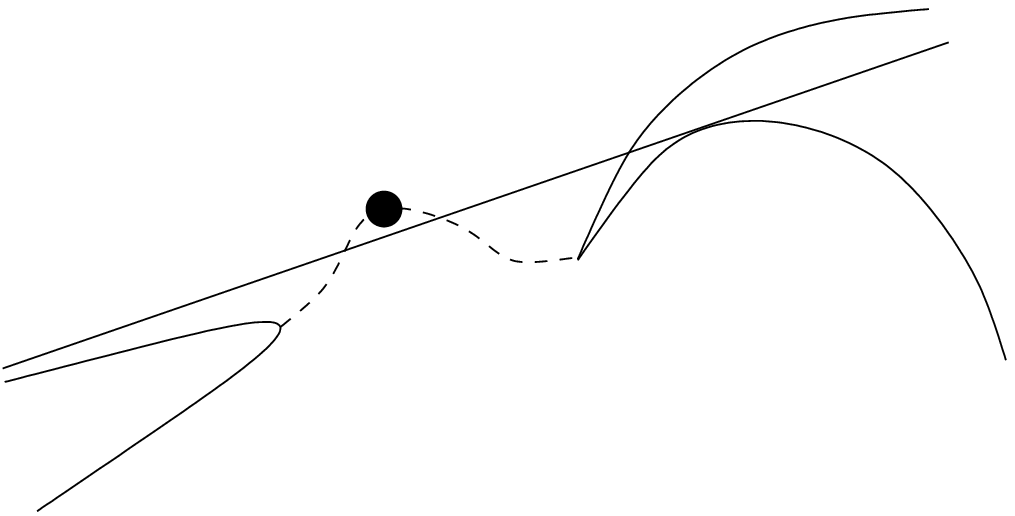}
\caption{Real part of $\cc_{\sqrt{2}}$}
\label{curvamas}
\end{figure}

\begin{figure}[ht]
\includegraphics[height=3cm]{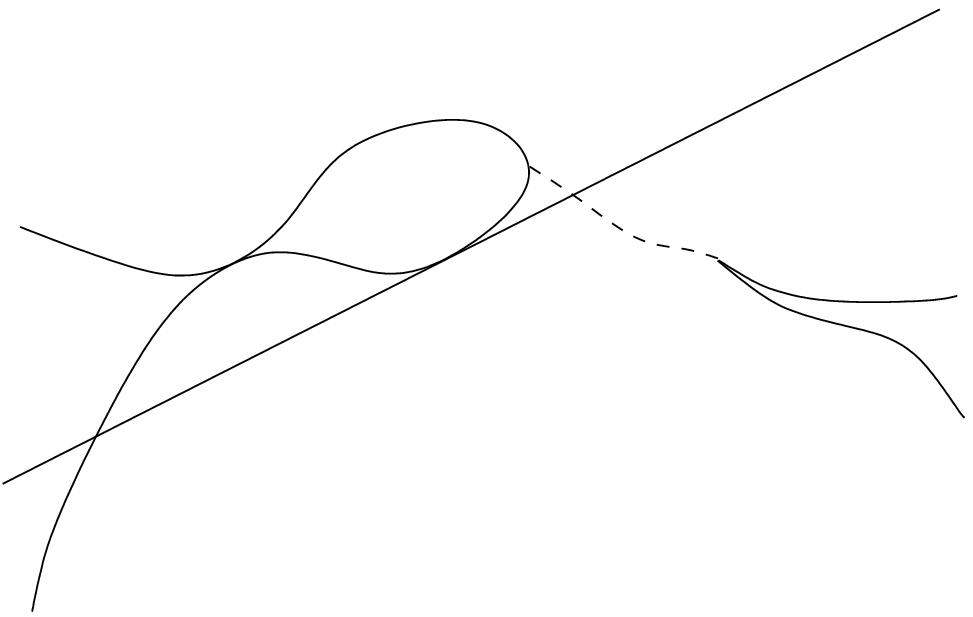}
\caption{Real part of $\cc_{-\sqrt{2}}$}
\label{curvamenos}
\end{figure}

Figures \ref{curvamas} and \ref{curvamenos}
show the real parts of $\cc_{\sqrt{2}}$
and $\cc_{-\sqrt{2}}$; we have drawn their
topological behavior. The dotted curves represent the real
parts of the imaginary solutions and the thick point
is the tacnode. The branch at infinity corresponding
to the $\be_6$ point is represented in both cases by
the branches of the curves going to $-\infty$.

Braid monodromy of $\cc_{\sqrt{2}}$ is computed from
figure \ref{curvamas}:
\begin{align*}
\alpha_7^{\sqrt{2}}&\mapsto\sigma_2^8\\
\alpha_1^{\sqrt{2}}&\mapsto\sigma_2^{4}*\sigma_1^2\\
\alpha_2^{\sqrt{2}}&\mapsto\sigma_2^{4}\sigma_1
*\sigma_2^3=\sigma_2^3*\sigma_1^3\\
\alpha_3^{\sqrt{2}}&\mapsto\sigma_2^{4}\sigma_1
\sigma_2\sigma_1^{-1}\sigma_2*\sigma_2^4=
\sigma_2*\sigma_1^4\\
\alpha_0^{\sqrt{2}}&\mapsto\sigma_2^{4}\sigma_1
\sigma_2\sigma_1^{-1}\sigma_2\sigma_1^2
\sigma_2^{-1}\sigma_1*\sigma_2=\sigma_1^{-3}*\sigma_2.
\end{align*}
Braid monodromy of $\cc_{-\sqrt{2}}$ is computed from
figure \ref{curvamenos}:
\begin{align*}
\alpha_2^{-\sqrt{2}}&\mapsto\sigma_2^3\\
\alpha_0^{-\sqrt{2}}&\mapsto\sigma_2\sigma_1^{-1}
\sigma_2*\sigma_1\\
\alpha_7^{-\sqrt{2}}&\mapsto\sigma_2\sigma_1^{-1}
\sigma_2\sigma_1*\sigma_2^8=
\sigma_2*\sigma_1^8\\
\alpha_3^{-\sqrt{2}}&\mapsto\sigma_2\sigma_1^{-1}
\sigma_2\sigma_1\sigma_2^4*\sigma_1^4=
\sigma_1^{-2}*\sigma_2^4\\
\alpha_1^{-\sqrt{2}}&\mapsto\sigma_2\sigma_1^{-1}
\sigma_2\sigma_1
\sigma_2^4\sigma_1^2*\sigma_2^2=\sigma_1^{-3}
*\sigma_2^2.
\end{align*}

Applying Zariski-Van Kampen Theorem \ref{zvk1}, we obtain
the fundamental group of the complement to the projective curve
$\cc_\beta\cup T_\beta$, where $T_\beta$ is the tangent line at $\be_6$.
In order to obtain $\pi_1(\bp^2\setminus \cc_\beta)$ we must
\emph{kill} a meridian about the line at infinity. 
Let $g_1,g_2,g_3$ denote the generators 
of the lexicographic geometric basis in a generic line
$z=K$, $K\gg 0$. Let $e$ denote the meridian about 
the point at infinity on this line, obtained by taking
$e:=(g_3\cdot g_2\cdot g_1)^{-1}$. On a neighbourhood 
of the center of the projection one can check that
\bdm
\left((e^{-1} g_3 e)\cdot g_3\cdot (e g_3 e^{-1})\cdot
g_3\cdot g_2\cdot g_1\right)^{-1}
\edm
is a meridian of the line at infinity.
Using GAP4 we obtain that both fundamental groups are 
isomorphic to $\bz \times SL(2;\bz/7\bz)$. 
Moreover, there is an isomorphism preserving meridians
--\,the image of a meridian in $\cc_{\sqrt{2}}$ is a 
meridian in $\cc_{-\sqrt{2}}$.
\smallbreak

\begin{thm}\label{noequiv} 
The braid monodromies of $\cc_{\sqrt{2}}$ and $\cc_{-\sqrt{2}}$
with respect to $\pi$ are not equivalent.
\end{thm}

\begin{proof}
It is enough to apply the final method described in strategy \ref{pasos}
using the representation
\bdm
B_3\to SL(2;\bz)\to SL(2;\bz/32\bz).
\edm
The orbits of both braid monodromies have size 15360 and
are disjoint as shown by the GAP4 program in the Appendix.
\end{proof}

Finally, we are in the position to prove theorem \ref{conj}.

\begin{proof}[Proof of Theorem \ref{conj}]
Let $\cc_1$ and $\cc_2$ in the statement of theorem \ref{conj}
correspond to $\cc_{\sqrt{2}}$ and $\cc_{-\sqrt{2}}$ respectively.
By theorems \ref{noequiv} and \ref{top} there is no 
homeomorphism $(\bp^2,\cc_1^\varphi)\to(\bp^2,\cc_2^\varphi)$
preserving orientations on $\bp^2$, $\cc_1$ and $\cc_2$.

Since both curves have real equations, they are invariant
by complex conjugation which preserves orientations
on $\bp^2$ but exchanges orientations on the curves.
Hence there is no homeomorphism
$(\bp^2,\cc_1^\varphi)\to(\bp^2,\cc_2^\varphi)$
preserving the orientation of $\bp^2$ but 
reversing it on $\cc_1$ and $\cc_2$.

Since algebraic braids always turn in the same direction
--\,that is, they are positive\,-- it is not possible to 
have a homeomorphism
$(\bp^2,\cc_1^\varphi)\to(\bp^2,\cc_2^\varphi)$
preserving the orientation of $\bp^2$ and 
only some, but not all, of the components of the curves.
\end{proof}

\appendix{\scshape Appendix}

In this appendix we provide the program source 
in GAP4 which produces the result stated in theorem \ref{noequiv}.
The authors include it for completeness and the text file 
can be distributed upon request. The execution of the program
took about ten hours on a Pentium III 866Mhz running with
GNU/Linux and sharing CPU time with other computations.

Here is the program:
\smallbreak
\noindent\# Some function definitions. Right conjugation.
\smallbreak

{\tiny\ttfamily
\begin{verbatimtab}[2]
cnj:=function(u,v)
        return u*v/u;
end;
\end{verbatimtab}
}
\noindent\# Simultaneous conjugation of a list.
{\tiny\ttfamily
\begin{verbatimtab}[2]
conjorbita:=function(lista,u)
        return List(lista,x->[x[1],x[2]^u]);
end;
\end{verbatimtab}
}

\smallbreak
\noindent\# Reverse product of a list.
\smallbreak

{\tiny\ttfamily
\begin{verbatimtab}[2]
prdct:=function(lista)
        local j,producto,n;
        n:=Length(lista);
        producto:=();
        for j in [1..n] do
                producto:=producto*lista[n-j+1][2];
        od;
        return producto;
end;
\end{verbatimtab}
}

\smallbreak
\noindent\# Action of  $\sigma_1$.
\smallbreak

{\tiny\ttfamily
\begin{verbatimtab}[2]
q1:=function(lista)
        local resultado,i,n;
        n:=Length(lista);
        resultado:=[];
        resultado[1]:=lista[2];
        resultado[2]:=[lista[1][1],cnj(lista[2][2],lista[1][2])];
        for i in [3..n] do
                resultado[i]:=lista[i];
        od;
        return resultado;
end;
\end{verbatimtab}
}

\smallbreak

\noindent\# Action of  $\sigma_1\cdot\ldots\cdot\sigma_{n-1}$.
\smallbreak
{\tiny\ttfamily
\begin{verbatimtab}[2]

q2:=function(lista)
        local resultado,i,n;
        n:=Length(lista);
        resultado:=[];
        resultado[1]:=lista[n];
        for i in [2..n] do
                resultado[i]:=[lista[i-1][1],cnj(lista[n][2],lista[i-1][2])];
        od;
        return resultado;
end;
\end{verbatimtab}
}

\smallbreak
\noindent\# This part produces the images of the standard
generators of the braid group by the given representation.
This part should be replaced for different representations.
\smallbreak

{\tiny\ttfamily
\begin{verbatimtab}[2]
m:=32;
R:=ZmodnZ(m);
fam:=ElementsFamily(FamilyObj(R));;
u:= ZmodnZObj(fam,1); 
um:= ZmodnZObj(fam,m-1); 
u0:= ZmodnZObj(fam,0); 
A:=[[u,u0],[um,u]];
B:=[[u,u],[u0,u]];
g1:=Group(A,B);
iso:=IsomorphismPermGroup(g1);
a:=Image(iso,A);
b:=Image(iso,B);
g:=Group(a,b);
\end{verbatimtab}
}

\smallbreak
\noindent\# This part describes both braid monodromies. 
This part should be replaced for different braid monodromies.
Note that the elements in the list have two entries. The first
one is a label to control the local topological singularity types.
{\tiny\ttfamily
\smallbreak

\begin{verbatimtab}[2]
pr:=[[1,b^8],[2,cnj(b^4,a^2)],[3,cnj(b^3,a^3)],
[4,cnj(b,a^4)],[5,b^(a^3)]];
otro:=[[3,b^3],[5,cnj(b/a*b,a)],[1,cnj(b,a^8)],
[4,(b^4)^(a^2)],[2,(b^2)^(a^3)]];
\end{verbatimtab}
}

\smallbreak
\noindent\# The reversed product is applied to obtain
pseudo Coxeter elements.
\smallbreak

{\tiny\ttfamily
\begin{verbatimtab}[2]

totpr:=prdct(pr);
tototro:=prdct(otro);
\end{verbatimtab}
}

\smallbreak
\noindent\# Verification to check if both pseudo Coxeter 
elements are conjugate. Then, the list of all conjugates 
to the second braid monodromy having the same pseudo 
Coxeter elements as the first one is produced.
\smallbreak
{\tiny\ttfamily
\begin{verbatimtab}[2]
Print(IsConjugate(g,tototro,totpr),"\n");
vale:=RepresentativeAction(g,tototro,totpr);
conjugar:=List(Elements(Centralizer(g,tototro)),x->x*vale);;
segundo:=Unique(List(conjugar,u->conjorbita(otro,u)));;
Sort(segundo);;
\end{verbatimtab}
}
\smallbreak
\noindent\# 
This part inductively constructs a subset, say $A$, containing 
the first braid monodromy up to conjugation and
stable by the action of the function $q_1$. Note that
it is enough to produce conjugations which preserve the
pseudo Coxeter element of the first braid monodromy.
Then it considers a subset $B$ of $A$ (up to the parameter $s$)
such that its image by the function $q_2$ is contained in $A$.
It chooses an element $x\in A\setminus B$ and applies 
$q_2$ to it. If its image is already in $A$ then $x$ is added to 
$B$. If $q_2(x)\notin A$ then its orbit by $q_1$ is added to $A$.
The program stops when $A=B$. In fact, at each step it checks for
common elements with the second braid monodromy and stops if 
there is any common element.
 \smallbreak
{\tiny\ttfamily
\begin{verbatimtab}[2]

cnm:=Centralizer(g,totpr);
lcnm:=Elements(cnm);
micnj:=function(el)
        return Unique(List(lcnm,x->conjorbita(el,x)));
end;
orbita:=[pr];
quedan:=ShallowCopy(orbita);
Sort(quedan);
s:=0;
t:=Length(orbita);
r:=Length(quedan);

while s<Length(orbita) do
        elemento:=ShallowCopy(orbita[t]);
        elmcnj:=micnj(elemento);
        elt:=ShallowCopy(orbita[t]);
        control:=not (elt in segundo);
        if not control then
                s:=Length(orbita);
                Print("Orbits are equal\n");
        fi;
        control0:=control;
        while control do
                elt1:=q1(elt);
                control:=not (elt1 in elmcnj);
                if control then
                        control0:=not (elt1 in segundo);
                        if control0 then
                                Add(orbita,elt1);
                                Add(quedan,elt1);
                                elt:=ShallowCopy(elt1);
                        else
                                s:=Length(orbita);
                                Print("Orbits are equal\n");
                                control:=false;
                        fi;
                fi;
        od;
        if control0 then
                t:=Length(orbita)+1;
                s:=s+1;
                Print(s," a ",t,"\n");
                elt:=ShallowCopy(orbita[s]);
                control:=true;
                Sort(quedan);
                while control do
                        elt1:=q2(elt);
                        control1:=not (elt1 in segundo);
                        control0a:=control1;
                        if not control1 then
                                s:=Length(orbita);
                                Print("Orbits are equal\n");
                                control:=false;
                        else
                                cntrlcnj:=true;
                                j:=1;
                                orbcnj:=micnj(elt1);
                                while cntrlcnj do
                                        elt1j:=orbcnj[j];
                                        control1:=not (elt1j in quedan);
                                        j:=j+1;
                                        cntrlcnj:=control1 and (not j>Length(orbcnj));
                                od;
                        fi;
                        if control0a then
                                if control1 then
                                        control:=false;
                                        Add(orbita,elt1);
                                        Add(quedan,elt1);
                                        Sort(quedan);
                                elif s<Length(orbita) then
                                        s:=s+1;
                                        Print(s,"b,\n");
                                        elt:=ShallowCopy(orbita[s]);
                                        RemoveSet(quedan,elt1);
                                else
                                        control:=false;
                                        s:=s+1;
                                        Print("Orbits are different\n");
                                fi;
                        fi;
                od;
        fi;
od;
\end{verbatimtab}
}

%\bibliographystyle{amsplain}
%\bibliography{biblio-ea}

\begin{thebibliography}{10}

\bibitem{ab:74}
H.~Abelson, \emph{Topologically distinct conjugate varieties with finite
  fundamental group}, Topology \textbf{13} (1974), 161--176.

\bibitem{acc:00}
E.~Artal~Bartolo, J.~Carmona, and J.I. Cogolludo, \emph{On sextic curves with
  big {M}ilnor number}, Preprint, November 2000.

\bibitem{acct:00}
E.~Artal~Bartolo, J.~Carmona, J.I. Cogolludo, and H.~Tokunaga, \emph{On curves
  with singular points in special position}, available at
  \texttt{arXiv:math.AG/0007152}, to appear in J. Knot Theory Ramifications.

\bibitem{art:47}
E.~Artin, \emph{Theory of braids}, Ann. of Math. (2) \textbf{48} (1947),
  101--126.

\bibitem{ass:78}
J.~Assion, \emph{Einige endliche {F}aktorgruppen der {Z}opfgruppen}, Math. Z.
  \textbf{163} (1978), no.~3, 291--302.

\bibitem{bir:74}
J.~S. Birman, \emph{Braids, links, and mapping class groups}, Princeton
  University Press, Princeton, N.J., 1974, Annals of Mathematics Studies, No.
  82.

\bibitem{brwj:86}
J.~S. Birman and B.~Wajnryb, \emph{Markov classes in certain finite quotients
  of {A}rtin's braid group}, Israel J. Math. \textbf{56} (1986), no.~2,
  160--178.

\bibitem{brk:88}
E.~Brieskorn, \emph{Automorphic sets and braids and singularities}, Braids
  (Santa Cruz, CA, 1986), Amer. Math. Soc., Providence, RI, 1988, pp.~45--115.

\bibitem{car:xx}
J.~Carmona, \emph{thesis}, preprint.

\bibitem{deg:90}
A.~I. Degtyar{\"e}v, \emph{Isotopic classification of complex plane projective
  curves of degree $5$}, Leningrad Math. J. \textbf{1} (1990), no.~4, 881--904.

\bibitem{GAP4}
The GAP~Group, Aachen, St~Andrews, \emph{{GAP -- Groups, Algorithms, and
  Programming, Version 4.2}}, 2000,
  \verb+(http://www-gap.dcs.st-and.ac.uk/~gap)+.

\bibitem{vk:33}
E.R.~van Kampen, \emph{On the fundamental group of an algebraic curve}, Amer.
  J. Math. \textbf{55} (1933), 255--260.

\bibitem{khku:01}
V.~Kharlamov and Vik.~S. Kulikov, \emph{Diffeomorphisms, isotopies, and braid
  monodromy factorizations of plane cuspidal curves}, available at
  \texttt{arXiv:math.AG/0104021}.

\bibitem{kl:88}
P.~Kluitmann, \emph{Hurwitz action and finite quotients of braid groups},
  Braids (Santa Cruz, CA, 1986), Amer. Math. Soc., Providence, RI, 1988,
  pp.~299--325.

\bibitem{kt:00}
Vik.~S. Kulikov and M.~Teicher, \emph{Braid monodromy factorizations and
  diffeomorphism types}, Izv. Ross. Akad. Nauk Ser. Mat. \textbf{64} (2000),
  no.~2, 89--120.

\bibitem{li:86}
A.~Libgober, \emph{On the homotopy type of the complement to plane algebraic
  curves}, J. Reine Angew. Math. \textbf{367} (1986), 103--114.

\bibitem{li:89}
\bysame, \emph{Invariants of plane algebraic curves via representations of the
  braid groups}, Invent. Math. \textbf{95} (1989), no.~1, 25--30.

\bibitem{li:98}
\bysame, \emph{Characteristic varieties of algebraic curves}, Preprint
  available at \texttt{arXiv:math.AG/ 9801070}, 1998.

\bibitem{mz:81}
B.~G. Moishezon, \emph{Stable branch curves and braid monodromies}, L.N.M. 862,
  Algebraic geometry (Chicago, Ill., 1980), Springer, Berlin, 1981,
  pp.~107--192.

\bibitem{se:64}
J.~P. Serre, \emph{Exemples de vari\'et\'es projectives conjugu\'ees non
  hom\'eomorphes}, C. R. Acad. Sci. Paris \textbf{258} (1964), 4194--4196.

\bibitem{waj:88}
B.~Wajnryb, \emph{Markov classes in certain finite symplectic representations
  of braid groups}, Braids (Santa Cruz, CA, 1986), Amer. Math. Soc.,
  Providence, RI, 1988, pp.~687--695.

\bibitem{waj:91}
\bysame, \emph{A braidlike presentation of ${\rm {s}p}(n,p)$}, Israel J. Math.
  \textbf{76} (1991), no.~3, 265--288.

\bibitem{zr:29}
O.~Zariski, \emph{On the problem of existence of algebraic functions of two
  variables possessing a given branch curve}, Amer. J. Math. \textbf{51}
  (1929), 305--328.

\end{thebibliography}
\providecommand{\bysame}{\leavevmode\hbox to3em{\hrulefill}\thinspace}

\end{document}